\documentclass[a4paper,12pt,leqno]{article}
\usepackage{latexsym}
\usepackage[all]{xy}

\usepackage{amssymb} 
\usepackage{amsmath} 
\usepackage{theorem}

\def\P{{\mathbf{P}}}
\def\Z{{\mathbb{Z}}}
\def\R{{\mathbb{R}}}
\def\A{{\mathcal{A}}}
\def\F{{\mathcal{F}}}
\def\G{{\mathcal{G}}}

\DeclareMathOperator{\Der}{Der}

\numberwithin{equation}{section}

\newcommand{\owari}{\hfill$\square$}

\theoremstyle{break}
\newtheorem{theorem}{Theorem}[section]
\newtheorem{prop}[theorem]{Proposition}
\newtheorem{cor}[theorem]{Corollary}
\newtheorem{lemma}[theorem]{Lemma}
\newtheorem{define}[theorem]{Definition}

\newtheorem{example}[theorem]{Example}

\newcommand{\xgraphAvertex}[1][****]{
\xgraphAVertex #1
}
\newcommand{\xgraphAVertex}[4]{
\if#3o\put(0,0){\circle{4}}\fi%
\if#3*\put(0,0){\circle*{4}}\fi%
\if#3.\put(0,0){\circle*{2.4}}\fi%
\if#4o\put(30,0){\circle{4}}\fi%
\if#4*\put(30,0){\circle*{4}}\fi%
\if#4.\put(30,0){\circle*{2.4}}\fi%
\if#2o\put(0,30){\circle{4}}\fi%
\if#2*\put(0,30){\circle*{4}}\fi%
\if#2.\put(0,30){\circle*{2.4}}\fi%
\if#1o\put(30,30){\circle{4}}\fi%
\if#1*\put(30,30){\circle*{4}}\fi%
\if#1.\put(30,30){\circle*{2.4}}\fi%
}
\newcommand{\xgraphA}[6]{
\if#1+\put(30,30){\line(-1,0){30}}\fi
\if#1.\qbezier[7](30,30)(15,30)(0,30)\fi
\if#1-\put(30,31){\line(-1,0){30}}\put(30,29){\line(-1,0){30}}\fi
\if#2+\put(0,0){\line(1,1){30}}\fi 
\if#2.\qbezier[10](0,0)(15,15)(30,30)\fi 
\if#2-\put(-0.7,0.7){\line(1,1){30}}\put(0.7,-0.7){\line(1,1){30}}\fi 
\if#3+\put(30,30){\line(0,-1){30}}\fi
\if#3.\qbezier[7](30,0)(30,15)(30,30)\fi 
\if#3-\put(29,30){\line(0,-1){30}}\put(31,30){\line(0,-1){30}}\fi
\if#4+\put(0,0){\line(0,1){30}}\fi 
\if#4.\qbezier[7](0,0)(0,15)(0,30)\fi 
\if#4-\put(-1,0){\line(0,1){30}}\put(1,0){\line(0,1){30}}\fi 
\if#5+\put(0,30){\line(1,-1){30}}\fi 
\if#5.\qbezier[10](30,0)(15,15)(0,30)\fi 
\if#5-\put(-0.7,29.3){\line(1,-1){30}}\put(0.7,30.7){\line(1,-1){30}}\fi 
\if#6+\put(0,0){\line(1,0){30}}\fi 
\if#6.\qbezier[7](0,0)(15,0)(30,0)\fi 
\if#6-\put(0,1){\line(1,0){30}}\put(0,-1){\line(1,0){30}}\fi 
\xgraphAvertex}

\title{
A primitive derivation and
logarithmic differential
forms of Coxeter arrangements}
\author{Takuro Abe\thanks{
Department of Mathematics, Kyoto University, 
Kitashirakawa-Oiwake-cho, Sakyo-Ku, 
Kyoto 606-8502, Japan.
email:abetaku@math.kyoto-u.ac.jp.} and 
Hiroaki Terao\thanks{
Department of Mathematics, Hokkaido University, 
Kita-10, Nishi-8, Kita-Ku, 
Sapporo, Hokkaido 060-0810, Japan.
email:terao@math.sci.hokudai.ac.jp.}
}
\date{} 

\begin{document}

\maketitle

\begin{abstract}
Let $W$ be a finite irreducible real
reflection group, which is a Coxeter group.
We explicitly construct a 
basis for the module of differential $1$-forms
with logarithmic poles 
along the Coxeter arrangement
by using
a primitive derivation. 
As a consequence,
we extend the Hodge filtration, 
indexed by nonnegative
integers, into 
a filtration indexed by all integers.
This filtration coincides with the filtration 
by the order of poles.
The results are  translated into the
derivation case.
\end{abstract}
\setcounter{section}{0}

\section{Introduction and main results}
Let $V$ be a Euclidean space of dimension $\ell$.
Let $W$ be a finite irreducible
reflection group (a Coxeter group)
acting on $V$.
%
The 
{\bf
Coxeter arrangement} $\A = \A(W)$
corresponding to $W$
is 
the set of reflecting hyperplanes.
We use \cite{OT} as a general
reference for arrangements.
For each $H\in \A$, choose a linear form 
$\alpha_{H} \in V^{*} $
such that $H = \ker (\alpha_{H} )$.
Their product $Q := \prod_{H\in \A} \alpha_{H} $,
which lies in the symmetric algebra
$S := \mbox{Sym}(V^*)$, is a defining polynomial for $\A$.
Let
$F:=S_{(0)}$ be
the quotient field of $S$.
Let $\Omega_{S} $ and
 $\Omega_{F}  $ denote the $S$-module of regular 
$1$-forms on $V$ and the $F$-vector space of
rational $1$-forms on $V$ respectively.  
The action of $W$ on $V$
induces 
the canonical actions of $W$ on
$V^{*}, S, F, \Omega_{S}$
and
$\Omega_{F}$,  
which enable us to consider
their $W$-invariant parts.
Especially let $R = S^{W} $ denote the
invariant subring of $S$.
%

In \cite{Z}, Ziegler
introduced the $S$-module of {\bf logarithmic $1$-forms} 
with poles of order $m \,\,(m\in \Z_{\geq 0} )$ 
along $\A$ by
\begin{align*}
\Omega(\A, m)
:=
\{
\omega\in\Omega_{F} \mid
&Q^{m} \omega \text{~and~} 
(Q/\alpha_{H})^{m} 
\left({d\alpha_{H}\wedge
\omega}\right) \\
&\text{~are both regular for all~} H\in\A
\}.
\end{align*} 
Note
$\Omega(\A, 0)=\Omega_{S}$. 
Define the total module of logarithmic $1$-forms
by
\[
\Omega(\A, \infty)
:=
\bigcup_{m\geq 0} \Omega(\A, m).
\]
In this article we study 
the total module
$\Omega(\A, \infty)$
of logarithmic $1$-forms 
and its
$W$-invariant part 
$\Omega(\A, \infty)^{W} $
by introducing 
a geometrically-defined filtration
indexed by $\Z$.

Let
$P_1, \cdots, 
P_{\ell}
\in R$ 
be
algebraically independent homogeneous polynomials 
with $\deg P_{1} \leq\dots\leq \deg P_{\ell}$,
which 
are called {\bf basic invariants},
 such that $R={\mathbb R}[P_1, \cdots, P_{\ell}]$
 \cite[V.5.3, Theorem 3]{bou}.
Define the {\bf primitive derivation} 
$D := \partial/\partial P_{\ell} : F \rightarrow F$.
%
Let $T := \{f \in R \mid Df = 0\}
=
\R [P_{1}, P_{2}, \dots , P_{\ell-1}].
$
Consider the $T$-linear connection (covariant derivative)
$$\nabla_{D}   : \Omega_{F}  \to 
\Omega_{F}$$ 
characterized by 
$
\nabla_{D} (f\omega) =
(Df) \omega +
f
(\nabla_{D} \omega)
\,\,\,
(f\in F, \omega\in \Omega_{F})$ 
and
$
\nabla_{D} (d\alpha) = 0
\,\,\,
(\alpha\in V^{*})$. 

In Section 2,
using the primitive derivation
 $D$, we 
explicitly construct 
logarithmic $1$-forms 
$$
\omega^{(m)}_{1},
\omega^{(m)}_{2},
\dots
,
\omega^{(m)}_{\ell}
$$
for each $m\in\Z$
satisfying
$\nabla_{D}\,\omega^{(2k+1)}_{j} = \omega^{(2k-1)}_{j}  
\ \ (k\in\Z, 1\le j\le\ell). $ 
The $1$-forms 
$\omega_1^{(m)},\ldots,\omega_\ell^{(m)}$
form a basis for the $S$-module
$\Omega(\A, -m)$
when $m\le 0$.
Thus it is natural to define
$\Omega(\A, -m)$ to be the $S$-module 
spanned by $
\{
\omega^{(m)}_{1},
\omega^{(m)}_{2},
\dots
,
\omega^{(m)}_{\ell}
\}
$
for all $m\in\Z.$ 
Let
$
\mathcal B_{k} :=\{
\omega^{(2k+1)}_{1},
\omega^{(2k+1)}_{2},
\dots
,
\omega^{(2k+1)}_{\ell}\}
$
for $k\in\Z$. 
The following two main theorems will be proved in Section 2:

\begin{theorem}
\label{0.4}

\hspace{5mm} 
(1)
The $R$-module 
$\Omega(\A, 2k-1)^{W} $
is free
with a basis
$\mathcal B_{-k} $ 
for $k\in\Z$.

(2)
The $T$-module 
$\Omega(\A, 2k-1)^{W} $
is free
with a basis
$\bigcup_{p\ge -k} \mathcal B_{p} $
for $k\in\Z$.

(3)
$
\mathcal B
:=
\bigcup_{k\in\Z} \mathcal B_{k} 
$ 
is a basis for 
$\Omega(\A, \infty)^{W} $
as a $T$-module.
\end{theorem}

\begin{theorem}
\label{0.1}

\hspace{5mm} (1) The $\nabla_{D} $ induces a $T$-linear automorphism  
$\nabla_{D} : 
\Omega (\A, \infty)^{W} 
\stackrel{\sim}{\rightarrow} 
\Omega (\A, \infty)^{W}. 
$ 

(2)
Define
$
\mathcal F_{0} := \bigoplus_{j=1}^{\ell} T \left(dP_{j}\right), 
\mathcal F_{-k} := \nabla_{D}^{k} \mathcal F_{0}$
and
$
\mathcal F_{k} := (\nabla_{D}^{-1})^{k} \mathcal F_{0}
\ \ (k>0)$.
Then
$
\Omega (\A, \infty)^{W} 
=
{
\bigoplus_{k\in\Z}}\, {\mathcal F}_{k}.
$

(3)
$
\Omega(\A, 2k-1)^{W}
=
\mathcal J^{(-k)}
$,
where
$
\mathcal J^{(-k)}
:= 
{\bigoplus_{p\geq -k}}\, \mathcal F_{p}
$ 
for $k\in\Z$.
\end{theorem}

\medskip

Let us briefly discuss our results
in connection with earlier researches.
Let
$\Der_{F} $ denote
the $F$-vector space
of $\R$-linear derivations of $F$ to itself. 
It is dual to $\Omega_{F} $. 
The inner product $I:V \times V \rightarrow \R$ 
induces
$I^*:V^* \times V^* \rightarrow \R$,
which
is canonically extended to a nondegenerate
$F$-bilinear form
$I^{*}  :
\Omega_{F} \times 
\Omega_F \rightarrow F$.
Define
an $F$-linear isomorphism 
\[
I^{*} : \Omega_{F} \rightarrow \Der_{F} 
\]
by
$I^{*} (\omega) (f) := I^{*} (\omega, df)
\,\,
(f\in F)$. 
Let 
$\mathcal G_{k} := I^{*} (\mathcal F_{k-1} )$ 
and
${\mathcal H}^{(k)} := I^{*} (\mathcal J^{(k-1)} )$ 
for $k\in\Z$. 
Thanks to Theorem \ref{0.1},
we have commutative
diagrams 
$$\xymatrix@R1pc{
\cdots \ar[r]^{\nabla_D} & \F_1 \ar[r]^{\nabla_D} \ar[d]_{I^{*}} & \F_0 \ar[r]^{\nabla_D} \ar[d]_{I^{*}}  &  \F_{-1} \ar[r]^{\nabla_D} \ar[d]_{I^{*}} 
&\F_{-2} \ar[r]^{\nabla_D} \ar[d]_{I^{*}} 
&\F_{-3} \ar[r]^{\nabla_D} \ar[d]_{I^{*}} 
&\F_{-4} \ar[r]^{\nabla_D} \ar[d]_{I^{*}} 
& \cdots\\
\cdots \ar[r]^{\nabla_D} & \G_2 \ar[r]^{\nabla_D}  & \G_1 \ar[r]^{\nabla_D}   &  \G_0    \ar[r]^{\nabla_D} 
&\G_{-1}    \ar[r]^{\nabla_D} 
&\G_{-2}    \ar[r]^{\nabla_D} 
&\G_{-3}    \ar[r]^{\nabla_D} 
& \cdots ,
} $$
$$\xymatrix@R1pc{
\cdots \ar[r]^{\nabla_D} 
&\mathcal J^{(1)} \ar[r]^{\nabla_D} \ar[d]_{I^{*}} 
&\mathcal J^{(0)} \ar[r]^{\nabla_D} \ar[d]_{I^{*}}  
&\mathcal J^{(-1)} \ar[r]^{\nabla_D} \ar[d]_{I^{*}} 
&\mathcal J^{(-2)} \ar[r]^{\nabla_D} \ar[d]_{I^{*}} 
&\mathcal J^{(-3)} \ar[r]^{\nabla_D} \ar[d]_{I^{*}} 
&\mathcal J^{(-4)} \ar[r]^{\nabla_D} \ar[d]_{I^{*}} 
& \cdots\\
\cdots \ar[r]^{\nabla_D} 
&\mathcal H^{(2)} \ar[r]^{\nabla_D}  
&\mathcal H^{(1)} \ar[r]^{\nabla_D}  
&\mathcal H^{(0)} \ar[r]^{\nabla_D}  
&\mathcal H^{(-1)} \ar[r]^{\nabla_D}  
&\mathcal H^{(-2)} \ar[r]^{\nabla_D}  
&\mathcal H^{(-3)} \ar[r]^{\nabla_D}  
& \cdots .
} $$
in which every $\nabla_{D}$ is a
$T$-linear isomorphism.
The objects in the left halves of the diagrams
were introduced by K. Saito who called
the decomposition
$\Der_{R} = \bigoplus_{k\geq 0} \,\mathcal G_{k}  $ 
the {\bf Hodge decomposition} 
and
the filtration 
$\Der_{R} = \mathcal H^{(0)} \supset \mathcal H^{(1)} \supset \dots$ 
the {\bf Hodge filtration}
in his groundbreaking work
\cite{Sa3, Sa4}.
They
are the key to
define the flat structure on the orbit space
$V/W$.  The flat structure is also called
the Frobenius manifold structure from the view point
of topological field theory \cite{dub1}.   

Our main theorems \ref{0.4} and \ref{0.1} are naturally
translated by
 $I^{*}$ 
into the corresponding
results concerning the $\mathcal G_{k} $'s 
and the $\mathcal H^{(k)}$'s 
in Section \ref{section3}. 
So we extend the Hodge decomposition and 
Hodge filtration,
{\bf
indexed by nonnegative integers}, 
to the ones 
{\bf
indexed by all integers.}
The Hodge
filtration
$\Der_{R} = \mathcal H^{(0)} \supset \mathcal H^{(1)} \supset \dots$ 
was proved to be equal to
 the contact-order
filtration \cite{T6}.
On the other hand,
Theorem \ref{0.1} (3)
asserts that
the filtration
$
\dots
\supset
\mathcal J^{(-1)} \supset 
\mathcal J^{(0)} 
=\Omega_{R} 
$,
indexed by nonpositive integers,
coincides with
the {\bf pole-order filtration} of 
the $W$-invariant part 
$\Omega(\A, \infty)^{W}$
of the total module 
$\Omega(\A, \infty)$
of logarithmic $1$-forms. 
This direction of researches is related with 
a generalized multiplicity $\textbf{m}:\A \rightarrow \Z$ 
and the associated logarithmic module $D\Omega(\A,\textbf{m})$ introduced in 
\cite{A4}.

In Section \ref{section4}, we will give explicit 
relations
of our
bases to the bases 
obtained in \cite{T4}, \cite{Y0}  and \cite{AY2}.

\section{Construction of a basis for $\Omega(\A, \infty) $ }
\label{constructionofbases} 
Let $
x_1,\ldots,x_\ell
$ denote a basis for $V^*$ and 
$P_1,\ldots,P_\ell$ homogeneous basic invariants
with $\deg P_{1} \leq\dots\leq \deg P_{\ell}$
:
$S^{W} = R = \R[P_1,\ldots,P_\ell]$. 
Let 
$\mathbf x := \left[
x_1,\ldots,x_\ell
\right]$ 
and
$\mathbf P := \left[
P_1,\ldots,
P_\ell
\right]$ 
be the corresponding row vectors.
Define
$A:=[I^*(x_i,x_j)]_{1 \le i,j \le \ell}
\in 
\mbox{GL}_{\ell} (\R)$
and
$
G:=[I^*(dP_i,dP_j)]_{1 \le i,j \le \ell}
\in \mbox{M}_{\ell,\ell}(R).
$
Then 
$
G=
J(\mathbf{P})^T A J(\mathbf{P}), 
$
where 
$
J(\mathbf{P})
:=
\left[
\displaystyle \frac{\partial P_j}{\partial x_i}
\right]_{1 \le i,j \le \ell}
$
is the 
Jacobian matrix.
It is well-known (e.g., \cite[V.5.5, Prop. 6]{bou}) that
$\det J(\P) \dot= Q,$
where
$\dot{=}$ stands for the equality up to a nonzero constant multiple.
 Let $\Der_{R} $ be the $R$-module of 
$\R$-linear derivations of $R$ to itself:
$\Der_{R} = \oplus_{i=1}^{\ell} R \,\left(\partial/\partial P_{i}\right).$  
Recall
the primitive derivation
$D
=
\partial/\partial P_{\ell} 
\in \Der_R$
and $T
=\ker(D:R \rightarrow R) = \R[P_1,\ldots,P_{\ell-1}]$. 
We will use the notation
$D[M]:=[D(m_{ij})]_{1 \le i,j \le \ell}$ 
for a matrix $M=[m_{ij}]_{1 \le i,j \le \ell} \in \mbox{M}_{\ell,\ell}(F)$.
The next Proposition is due to K. Saito
\cite[(5.1)]{Sa3} \cite[Corollary 4.1]{dub1}:

\begin{prop}
$
D[G]
\in
\mbox{GL}_{\ell} (T)
$,
that is,
$D^2[G]=0$ and $\det D[G] \in  \R^{\times} $.  
\label{DG}
\end{prop}

Now let us give a
key definition of this article, 
which generalizes the
matrices introduced in  
\cite[Lemma 3.3]{T4}. 

\begin{define}
The matrices $B=B^{(1)}$ and $B^{(k)}\ (k \in \Z)$ are defined by 
\begin{eqnarray*}
B:= J(\mathbf{P})^T A D[J(\mathbf{P})],\,\,\,\,
B^{(k)}:=kB+(k-1)B^T. 
\end{eqnarray*}
In particular, 
$
D[G]=B+B^T=B^{(k+1)}-B^{(k)}
$
for all $k \in \Z$. 
\label{matrixB}
\end{define}

\begin{lemma}
$
B^{(k)} \in \mbox{GL}_{\ell}(T)$ 
for all $k \in \Z$,
that is,
$D\left[B^{(k)}\right]=0$
and
 $\det B^{(k)}\in  \R^{\times} $.
\label{ngB}
\end{lemma}

\noindent
\textbf{Proof}. 
If $k\ge 1$, then the statement is proved 
in \cite[3.3 and 3.6]{T4} 
and \cite[Lemma 2]{T6}. Suppose $k \le 0$. 
Since
\begin{align*}
B^{(1-k)}=(1-k)B+(-k)B^T=-\{ kB + (k-1) B^T \}^T= - (B^{(k)})^T,
\end{align*}
we obtain
$B^{(k)} = - (B^{(1-k)} )^{T} \in \mbox{GL}_{\ell}(T) $ because $1-k\ge 1$. 
\owari

\medskip

The following Lemma is in
\cite[pp. 670, Lemma 3.4 (iii)]{T4}:

\begin{lemma} 
\hspace{5mm} (1)
$\det
J(D^k[\mathbf{x}])
\dot=
Q^{-2k},
$ 
where
$
J(D^k[\mathbf{x}]) := 
\left[\partial D^k (x_{j} )/\partial x_{i} \right]_{1\le i, j\le \ell}
\ (k\ge 1)$.


(2) $D[J(\P)]=-J(D[{\mathbf x}]) J(\P)$
and thus $\det D[J(\P)]\dot= Q^{-1}$. 
\label{detJDk}
\end{lemma}

\begin{define}
Define $\{R_{k}\}_{k \in \Z} \subset \mbox{M}_{\ell,\ell}(F)$ by 
\begin{eqnarray*}
R_{1-2k}:&=&D^k[J(\mathbf{P})]\ (k \ge 0),\\
R_{2k-1}:&=&(-1)^k J(D^k[\mathbf{x}])^{-1} D[J(\mathbf{P})]\ (k \ge 1),\\
R_{2k}:&=&(-1)^k J(D^k[\mathbf{x}])^{-1} \ (k \ge 0),\\
R_{-2k}:&=&D^{k+1}[J(\mathbf{P})] D[J(\mathbf{P})]^{-1} \ (k \ge 0).
\end{eqnarray*}
In particular, $R_{1} = J(\P)$,
 $R_{0} = I_{\ell} $ and $R_{-1}=D[J(\mathbf{P})]$.
\label{matrixR}
\end{define}

The following Proposition is fundamental.

\begin{prop}
For $k \in\Z$, we have

(1) 
$\det R_{k} \dot= Q^{k}$,

(2) 
$R_{2k} = R_{2k-1} D[J(\P)]^{-1} = R_{2k-1} B^{-1} J(\P)^T A,
$ 

(3) 
$R_{2k+1}  = 
R_{2k} J(\P) (B^{(k+1)})^{-1} B,
$

(4) 
$R_{2k+1}=R_{2k-1}B^{-1}G(B^{(k+1)})^{-1} B$,
and

(5)
$D[R_{2k+1} ] = R_{2k-1}. $ 
\label{inductive}
\end{prop}

\noindent
\textbf{Proof}. 
(2) is immediate from Definition \ref{matrixR}
because 
$B^{-1} J(\P)^{T}  A = D[J(\P)]^{-1} $. 

(4) 
Let $k\ge 1$. 
Recall  the original definition of $B^{(k)}$ in \cite[Lemma 3.3]{T4} given by 
\begin{eqnarray*}
B^{(k+1)}
=
-J(\mathbf{P})^T A J(D^{k+1}[\mathbf{x}]) J(D^k[\mathbf{x}])^{-1} J(\mathbf{P}).\end{eqnarray*}
Compute
\begin{eqnarray*}
R_{2k-1}^{-1}R_{2k+1}&=& 
-D[J(\mathbf{P})]^{-1} J(D^k[\mathbf{x}]) J(D^{k+1}[\mathbf{x}])^{-1} D[J(\mathbf{P})]\\
&=&-D[J(\mathbf{P})]^{-1} A^{-1} J(\mathbf{P})^{-T} J(\mathbf{P})^T A
J(\mathbf{P}) J(\mathbf{P})^{-1} \\
&\ & ~~~~~~J(D^k[\mathbf{x}]) J(D^{k+1}[\mathbf{x}])^{-1}
A^{-1} J(\mathbf{P})^{-T} J(\mathbf{P})^T A D[J(\mathbf{P})] \\
&=& B^{-1} G(B^{(k+1)})^{-1} B. 
\end{eqnarray*}
Next we will show 
that
$$
D^{k+1}[J(\mathbf{P})]=D^k[J(\mathbf{P})]B^{-1} B^{(1-k)} G^{-1} B
$$
for $k\ge 0$ by an
induction
on $k$. 
When $k=0$
we have
\begin{eqnarray*}
J(\mathbf{P})B^{-1} B^{(1)} 
G^{-1}B&=&J(\mathbf{P})J(\mathbf{P})^{-1}A^{-1}J(\mathbf{P})^{-T} 
J(\mathbf{P})^{T}AD[J(\mathbf{P})]
=D[J(\mathbf{P})].
\end{eqnarray*}
Next assume $k>0$. 
Compute
\begin{eqnarray*}
D^{k+1}[J(\mathbf{P})]&=&D[D^{k}[J(\mathbf{P})] ]
= D[D^{k-1}[J(\mathbf{P})]B^{-1}B^{(2-k)} G^{-1}B]\\
&=& D^{k}[J(\mathbf{P})]B^{-1}B^{(2-k)}G^{-1}B+
D^{k-1}[J(\mathbf{P})]B^{-1}B^{(2-k)} D[G^{-1}] B\\
&=& D^{k}[J(\mathbf{P})]B^{-1} \{B^{(2-k)} -D[G]\}G^{-1}B\\
&=&D^{k}[J(\mathbf{P})]B^{-1}B^{(1-k)}G^{-1}B,
\end{eqnarray*}
where, in the above, we used the induction hypothesis
$$
D^{k}[J(\mathbf{P})]=D^{k-1}[J(\mathbf{P})]B^{-1}B^{(2-k)}G^{-1}B,
$$
a general formula 
$$
D[G^{-1}]=-G^{-1}D[G]G^{-1}
$$
and 
$$
D[G]=B+B^T=B^{(2-k)}-B^{(1-k)}.
$$
This implies 
$R_{-2k-1} = R_{-2k+1} B^{-1} B^{(1-k)} G^{-1} B$ 
which proves (4).

(3) follows from (2) and (4) because $G = J(\P)^{T} A J(\P)$.

(1) Since
$
\det B^{(k)} \in\R^{\times},
$
$\det J(D^{k}[\mathbf x] )\dot= Q^{-2k}$
and
$
\det D[J(\P)] \dot= Q^{-1}$
by Lemma \ref{ngB} and
Lemma \ref{detJDk}, 
(1) is proved.

(5) 
follows from the following computation:
\begin{eqnarray*}
D[R_{2k+1}]B^{-1}&=&D[R_{2k+1}B^{-1}]
=D[R_{2k-1}B^{-1}G(B^{(k+1)})^{-1}]\\
&=&\{D[R_{2k-1}]B^{-1}G
+R_{2k-1}B^{-1}D[G]\}
(B^{(k+1)})^{-1}\\
&=&\{R_{2k-3} B^{-1}G+R_{2k-1} B^{-1} (B^{(k+1)}-B^{(k)})\}(B^{(k+1)})^{-1}\\
&=&\{R_{2k-1}B^{-1}B^{(k)}+R_{2k-1} B^{-1} (B^{(k+1)}-B^{(k)})\}(B^{(k+1)})^{-1}\\
&=&R_{2k-1}B^{-1}. \hspace{78mm} \square
\end{eqnarray*}

\begin{define}
For $m\in\Z$ define
$
\omega_1^{(m)},\ldots,\omega_\ell^{(m)}
\in
\Omega_{F}$
by 
$$
[
\omega_1^{(m)},\ldots,\omega_\ell^{(m)}
]
:=[dx_1,\ldots,dx_\ell]R_{m}.
$$
When $m = 2k+1\ \ (k\in\Z)$, let
\[
\mathcal B_{k}
:=
\{
\omega_1^{(2k+1)},\ldots,\omega_\ell^{(2k+1)}
\}. 
\]
\label{omega}
\end{define}

For example, $\omega^{(1)}_{j} = dP_{j}  $ for $1\le j\le \ell$
and $\mathcal B_{0} =\{dP_{1}, \dots , dP_{\ell} \}$ 
because 
$$
[\omega_1^{(1)},\ldots,\omega_\ell^{(1)}]
=[dx_1,\ldots,dx_\ell]J(\mathbf{P})=[dP_1,\ldots,dP_\ell].
$$


\begin{prop}
The subset
$$
\mathcal B
:= 
\bigcup_{k\in\Z} \mathcal B_{k} 
=
\{
\omega_{j}^{(2k+1)}
\mid
1 \le j \le \ell, \ k \in \Z
\}
$$
of 
$\Omega_{F}$ 
is linearly
independent over $T$.
\label{independent}
\end{prop}

\noindent
\textbf{Proof}. 
Assume 
$$
\sum_{k \in \Z} [\omega_1^{(2k+1)},\ldots,\omega_\ell^{(2k+1)}]{\mathbf g}^{(2k+1)}=0
$$
with
${\mathbf g}^{(2k+1)} =[g^{(2k+1)}_1,\ldots
g^{(2k+1)}_\ell]^T \in T^\ell,\ k \in \Z$
 such that 
there exist
integers $d$ and $ e$ such that $d \ge e$, 
${\mathbf g}^{(2d+1)}\neq 0,{\mathbf g}^{(2e+1)} \neq 0$ and 
${\mathbf g}^{(2k+1)}=0$ for all $k >d$ and $k <e$. 
Then 
$$
0=\sum_{k=e}^d [dx_1,\ldots,dx_\ell] R_{2k+1}{\mathbf g}^{(2k+1)}
$$
implies that 
$$
0=\sum_{k=e}^d R_{2k+1}{\mathbf g}^{(2k+1)}.
$$
By Proposition \ref{inductive} (4), there exist $(\ell \times \ell)$-matrices 
$H_{2k+1}\ (e \le k \le d)$ such that 
$$
R_{2k+1}=R_{2e+1} H_{2k+1}\ (e \le k \le d)
$$
and
$H_{2k+1}$ can be expressed as a product of 
$(k-e)$
copies of $G$ and 
matrices belonging to $\mbox{GL}_{\ell}(T)$. 
Since $\det (R_{2e+1}) \neq 0$
by Proposition \ref{inductive}  (1),  
$$
0=\sum_{k=e}^d H_{2k+1} {\mathbf g}^{(2k+1)}.
$$
Note
$D^{d-e}[H_{2k+1} ] = 0
\ (k < d)$ by Proposition
\ref{DG} and Lemma \ref{ngB}. 
Applying $D^{d-e}$ to the above, we thus
obtain
$$
D^{d-e}[H_{2d+1}]{\mathbf g}^{(2d+1)}=0.
$$
Since the matrix
$D^{d-e}[H_{2d+1}]$,
which is 
a product of 
$(d-e)$ copies of $D[G]$ and 
matrices in $\mbox{GL}_{\ell}(T)$,
is
nondegenerate,
we get
${\mathbf g}^{(2d+1)}=0$, which is a contradiction. \owari

\medskip

\begin{prop}
$\nabla_D \, \omega_j^{(2k+1)}=\omega_j^{(2k-1)}
\ (k \in \Z,\ 1\le j\le\ell)$.
%
\label{omegainductive}
\end{prop}

\noindent
{\bf Proof.}
By Proposition \ref{inductive} (5)
we have
\begin{align*} 
&~~~\left[
\nabla_{D} \,\omega_{1}^{(2k+1)},
\dots,
\nabla_{D} \,\omega_{\ell}^{(2k+1)}
\right]
=
\left[
dx_{1}, \dots, dx_{\ell}
\right]
D[R_{2k+1}]\\ 
&=
\left[
dx_{1}, \dots, dx_{\ell}
\right]
R_{2k-1} 
=
\left[
\omega_{1}^{(2k-1)},
\dots,
\omega_{\ell}^{(2k-1)}
\right].\ \ \ \square
\end{align*}

Recall
\begin{eqnarray*} 
\Omega(\A, \infty)
:&=&
\bigcup_{m\geq 0} \Omega(\A, m)\\
&=&
\{
\omega\in \Omega_{F} \mid
Q^{m} \omega \in \Omega_{S} \text{~ for some $m>0$ and~}\\
&~&~~~~~~~~~~~~~~d\alpha_{H} \wedge \omega 
\text{~is regular at generic points on~}H\\
&~&~~~~~~~~~~~~~\text{~for each~} H\in\A 
\}.
\end{eqnarray*}

\begin{lemma}
$\nabla_{D} (\Omega(\A, m)^{W} ) \subseteq \Omega(\A, m+2)^{W} $
for $m > 0$.
\label{closedinfty} 
\end{lemma}

%
%
%

\noindent
\textbf{Proof}. 
%
%
Choose $H\in\A$ arbitrarily and fix it.
Pick an orthonormal basis 
$\alpha_{H} = x_{1}, x_{2}, \dots , x_{\ell}   $ 
for $V^{*} $.
Let $s = s_{H} \in W$ be the orthogonal reflection
through $H$.  Then 
$s(x_{1} ) = - x_{1}, s(x_{i} )= x_{i} \,\,(i\geq 2) $,
$s(Q) = -Q$.
Let 
\[
\omega
=
\sum_{i=1}^{\ell} (f_{i}/Q^{m}) dx_{i} 
\in
\Omega(\A, m)^{W} 
\]
with each $f_{i} \in S$.
Then
\[
\nabla_{D} \,\omega
=
\sum_{i=1}^{\ell} D(f_{i}/Q^{m}) dx_{i} 
\]
is $W$-invariant with 
poles of order $m+2$ at most.
The $2$-form 
\[
(Q/x_{1})^{m} dx_{1} \wedge \omega
=
\sum_{i=2}^{\ell} (f_{i}/x_{1}^{m}) dx_{1} \wedge
dx_{i} 
\]
is regular because $\omega\in \Omega(\A, m)^{W}$.    
Let $i\ge 2$. 
Then $f_{i} \in x_{1}^{m} S$.
This implies that
 $g_{i} := Q^{m+2} D(f_{i}/Q^{m}) \in x_{1}^{m+1}  S$.
It is enough to show
$g_{i} \in x_{1}^{m+2} S $ because
\[
(Q/x_{1})^{m+2} 
dx_{1} \wedge\nabla_{D}\, \omega
=
\sum_{i=2}^{\ell} (g_{i}/x_{1}^{m+2}) dx_{1} \wedge dx_{i}. 
\]
When $m$ is odd, we have
$s(g_{i} )=s(Q^{m+2} D(f_{i}/Q^{m}))=-g_{i} $.
Thus
$g_{i} \in x_{1}^{m+2} S $.
When $m$ is even, we have
$s(g_{i} )=s(Q^{m+2} D(f_{i}/Q^{m}))=g_{i} $.
Thus
$g_{i} \in x_{1}^{m+2} S $.

%
\owari
\medskip

\begin{lemma}
$\mathcal B_{-k} \subset \Omega(\A, 2k-1)^{W} 
$ for $k \ge 1$. 
\label{omegasolution1}
\end{lemma}

\noindent
{\bf Proof.}
We will show by an induction on 
$k$.   
Fix $1\le j\le \ell$. 
Recall
$
\omega^{(-1)}_{j}
=
\nabla_{D}\, dP_{j}
$ 
by Proposition \ref{omegainductive}.
Since $dP_{j} \in\Omega(\A, 0)^{W} $,
we have
$\nabla_{D}\, dP_{j}\in\Omega(\A, 2)^{W} $
by Lemma \ref{closedinfty}.
On the other hand,  $\nabla_{D}\, dP_{j}$
has poles of order one at most because 
$dP_{j} $ is regular.  Thus
$\omega^{(-1)}_{j}\in\Omega(\A, 1)^{W} $. 
The induction proceeds by
Proposition \ref{omegainductive}
and Lemma
\ref{closedinfty}.
\owari

\medskip

We extend the definition of
$
\Omega(\A, m) 
$ to the case when $m$ is a negative integer:
\[
\Omega(\A, m) := \bigoplus_{j=1}^{\ell} S \, \omega^{(-m)}_{j}
\ \ \ (m<0).   
\]
%

\begin{theorem} 
$\Omega(\A, m)$ is a free $S$-module
with a basis 
$
\omega^{(-m)}_{1},
\omega^{(-m)}_{2},
\dots,
\omega^{(-m)}_{\ell}
$ for $m\in\Z.$   
\end{theorem}

\medskip

\noindent
{\bf Proof.}
 {\it Case 1.} When $m<0$ this is nothing but the definition.

 {\it Case 2.}
Let $m=2k-1$ with $k\ge 1$. 
Recall $\mathcal B_{-k} \subset \Omega(\A, 2k-1)^{W} $
from Lemma \ref{omegasolution1}
and $\det R_{1-2k} \dot= Q^{1-2k} $ by
Proposition \ref{inductive} (1).
Thus we have
\begin{eqnarray*} 
\omega_{1}^{(-2k+1)}
\wedge
\omega_{2}^{(-2k+1)}
\wedge
\dots
\wedge
\omega_{\ell}^{(-2k+1)}
&=&
\left(\det R_{1-2k}\right)
dx_{1} \wedge dx_{2} \wedge \dots \wedge dx_{\ell}\\ 
&\dot=&
Q^{1-2k} (dx_{1} \wedge dx_{2} \wedge \dots \wedge dx_{\ell}).
\end{eqnarray*} 
This shows that
$\mathcal B_{-k} $ is an $S$-basis
for $\Omega(\A, 2k-1)$
by Saito-Ziegler's criterion
\cite[Theorem 11]{Z}. 

 {\it Case 3.} Let $m=2k$ with $k\ge 0$. 
When $k=0$, the assertion is obvious because
$\omega^{(0)}_{j} = dx_{j}  $ and $\Omega(\A, 0) = \Omega_{S} $.
Let $k \ge 1$. 
By Proposition \ref{inductive} (2) we have
\begin{align*}
\left[
\omega^{(-2k)}_{1},
\dots,
\omega^{(-2k)}_{\ell}
\right]
&=
\left[
dx_{1},
\dots,
dx_{\ell}
\right]
R_{-2k} 
=
\left[
dx_{1},
\dots,
dx_{\ell}
\right]
R_{-2k-1}
B^{-1} J(\P)^T A\\
&=
\left[
\omega^{(-2k-1)}_{1},
\dots,
\omega^{(-2k-1)}_{\ell}
\right]
B^{-1} J(\P)^T A.
\end{align*} 
This implies that
$\omega^{(-2k)}_{1},
\dots,
\omega^{(-2k)}_{\ell}
$ 
lie in $\Omega(\A, 2k+1)$ by 
Lemma \ref{omegasolution1}.
By Proposition \ref{inductive} (3) we have
\begin{align*}
Q^{2k} R_{-2k} =
Q^{2k-1} 
R_{-2k+1} B^{-1} B^{(-k+1)} Q J(\P)^{-1}.
\end{align*} 
Since both 
$Q^{2k-1} 
R_{-2k+1} $ 
and
$Q J(\P)^{-1}$ belong to $M_{\ell, \ell} (S) $,
so does  $Q^{2k} R_{-2k}$.
In other words, the differential forms 
$\omega^{(-2k)}_{1},
\dots,
\omega^{(-2k)}_{\ell}
$ 
have poles of order at most $2k$ 
along $\A$. 
Since it is easy to see that
$\Omega(\A, 2k) = \Omega(\A, 2k+1) \cap (1/Q^{2k})\Omega_{S}$,
we know that
$\omega_j^{(-2k)}$ belongs to $\Omega(\A, 2k)$
for each $j$. 
We can apply Saito-Ziegler's criterion
\cite[Theorem 11]{Z}
to conclude
that
$\{\omega^{(-2k)}_{1},
\dots,
\omega^{(-2k)}_{\ell}\}
$ 
is a basis for $\Omega(\A, 2k)$ over $S$
because
$
\det R_{-2k} 
\dot=
Q^{-2k} 
$
by Proposition \ref{inductive} (1).
\owari

\medskip

We are now ready to prove
Theorems \ref{0.4} and \ref{0.1}. 

\medskip

\noindent
{\bf Proof of Theorem \ref{0.4}.}

(1) It is enough to show that
$\mathcal B_{-k} $ spans 
$\Omega(\A, 2k-1)^{W} $
over $R$.
Express an arbitrary element
$\omega\in\Omega(\A, 2k-1)^{W} $ 
as
\[
\omega
=
\sum_{j=1}^{\ell} f_{j} \omega^{(-2k+1)}_{j}   
\]
with each $f_{j} \in S$.
For any $s\in W$, get
\[
0 = \omega - s(\omega)
=
\sum_{j=1}^{\ell} \left[f_{j} - s(f_{j} )\right] \omega^{(-2k+1)}_{j}.   
\]
 Since $\mathcal B_{-k} $ is linearly independent
over $F$, we obtain $f_{j} \in S^{W} = R$.

(2)
Let $d_{j} := \deg P_{j} $ 
and
$m_{j} := d_{j} - 1$ for $1\le j\le \ell$. 
Let $h := d_{\ell} $ denote the Coxeter number.
Define the degree of a
homogeneous 
rational $1$-form by  
\[
\deg (\sum_{i=1}^{\ell}  f_{i} \,
dx_{i} ) = d 
\Longleftrightarrow
f_{i} = 0 \text{~or~} \deg f_{i} = d \ \ (1\le i \le\ell).
\]
Then 
\[
\deg \omega_{j}^{(2k+1)} = m_{j} + kh. 
\]
%
%
Recall that
$\mathcal B $ is linearly independent over $T$
by Proposition \ref{independent}.
Let $M_{-k} $ denote the free $T$-module 
 spanned
by  $\bigcup_{p\ge -k} \mathcal B_{p} $.
Recall that $\Omega(\A, 2k-1)^{W} $ is a free $R$-module with
a basis $\mathcal B_{-k}$
by  (1).
If $p\ge -k$, then  $R_{2p+1} = R_{-2k+1} H$ 
with a certain matrix $H\in M_{\ell, \ell} (R)$
because of Proposition \ref{inductive} (4).
This implies that $M_{-k} \subseteq 
\Omega(\A, 2k-1)^{W}
$.  
Use a Poincar\'e series argument to prove that
they are equal:
\begin{align*}
{\rm Poin}(M_{-k}, t)
&=
{
(1-t^{d_{1} } )^{-1} 
\dots
(1-t^{d_{\ell-1} } )^{-1} 
}  
\sum_{p\ge -k}
\left(
{
t^{m_{1} +ph} +
\dots
t^{m_{\ell} +ph} 
}
\right)\\
&=
{
(1-t^{d_{1} } )^{-1} 
\dots
(1-t^{d_{\ell} } )^{-1} 
}  
\left(
{
t^{m_{1} -kh} +
\dots
t^{m_{\ell} -kh} 
}
\right)\\
&=
{\rm Poin}(
\Omega(\A, 2k-1)^{W}, t).
\end{align*}
Therefore 
$M_{-k} = \Omega(\A, 2k-1)^{W}$.

  (3)
Thanks to Proposition \ref{independent},
it is enough to prove
that $\mathcal B$ spans
$\Omega(\A, \infty)^{W} $ over $T$.
Let $\omega\in\Omega(\A, \infty)$.
Then $\omega\in\Omega(\A, 2k-1)^{W} $
for some $k\ge 1$.
By (2) and (3) 
we conclude
that $\omega$ is a linear combination of
$\bigcup_{p\ge -k} \mathcal B_{p} $ with coefficients in $T$.
This shows that $\mathcal B$ spans $\Omega(\A, \infty)$
over $T$. 
\owari 

\medskip

\noindent
{\bf Proof of Theorem \ref{0.1} (1)}.
By Proposition \ref{omegainductive},
$$\nabla_{D} : \Omega(\A, \infty)^{W}
\rightarrow \Omega(\A, \infty)^{W}$$ induces a bijection 
$\nabla_{D} : \mathcal B \rightarrow \mathcal B.$ 
Apply Theorem \ref{0.4} (3) to prove that 
$\nabla_{D}$ is a $T$-isomorphism.
\owari

\medskip

Let $\nabla_{D}^{-1} : \Omega(\A, \infty)\rightarrow
 \Omega(\A, \infty)$ denote the inverse $T$-isomorphism.

\begin{define}
For $k \in \Z$,  define
$$
\F_{0} := \bigoplus_{j=1}^{\ell} T\,\left(dP_{j}\right), \ \ \
\F_{-k}:= \nabla_{D}^{k} (\F_{0} ) \ \ (k > 0), \ \ \ 
\F_{k}:= \left(\nabla_{D}^{-1}\right)^{k}  (\F_{0} ) \ \ (k > 0). 
$$ 
\label{defF}
\end{define}

Thus $\nabla_{D} $ induces a $T$-isomorphism 
$\nabla_{D}:
 \F_{k} \tilde{\rightarrow} \F_{k-1} $
for each $k\in \Z$.
Since
$\nabla_{D} $ induces a bijection 
$\nabla_{D} : \mathcal B_{k}  \rightarrow \mathcal B_{k-1}$ 
by Proposition \ref{omegainductive},
each $\F_k$ is a free $T$-module of rank $\ell$ with a basis
$
\mathcal B_{k} = \{
 \omega_j^{(2k+1)} \mid 1\le j\le\ell\}.
$

\medskip

\noindent
{\bf Proof of Theorem \ref{0.1} (2) and (3)}.

(2)
By Theorem \ref{0.4} (3), 
$
\mathcal B
=
\bigcup_{k\in\Z} \mathcal B_{k} 
$ 
is a basis for 
$\Omega(\A, \infty)^{W} $
as a $T$-module.
%
On the other hand,
each $\F_{k} $ has a basis 
$\mathcal B_{k} $ over $T$ for each 
$k\in\Z$. 

(3)
By Theorem \ref{0.4} (2),  
$\mathcal J^{(-k)} = \Omega(\A, 2k-1)^{W}$.
\owari

\medskip

\begin{example} 
Let $\A$ be the $B_{2}$ type
arrangement defined by 
$Q = xy(x+y)(x-y)$ corresponding 
to the Coxeter group of type $B_{2} $.
Then $P_{1} = (x^{2} + y^{2} )/2,\
P_{2} = (x^{4} + y^{4} )/4$ are basic invariants.
Then
$T = \R[P_{1} ]$ and $R = \R[P_{1}, P_{2}]$.
 Let
$$
\omega=
(x^{4} + y^{4})(\frac{dx}{x} + \frac{dy}{y}) \in \Omega(\A, 1)^{W}.
$$
The unique
decomposition of $\omega$ corresponding to
the decomposition 
$
\Omega(\A, 1)^{W} =
\mathcal J^{(-1)}
= 
\mathcal F_{-1}\oplus
\mathcal F_{0}\oplus
\mathcal F_{1}\oplus
\dots
$
is explicitly given by: 
\[
\omega = -8 P_{1}^{3} \omega^{(-1)}_{1} 
+(8/3) P_{1}^{2} \omega^{(-1)}_{2}    
-4 P_{1} \omega^{(1)}_{1}    
+2 \omega^{(1)}_{2}   
\in
\mathcal F_{-1} 
\oplus
\mathcal F_{0}   
\]
by an easy calculation.
\end{example}

\begin{cor} 
The
$\nabla_{D} : \Omega(\A, \infty)^{W} 
\rightarrow
\Omega(\A, \infty)^{W} $ 
induces 
an $T$-isomorphism 
$$
\nabla_{D} 
:
\Omega(\A, 2k-1)^{W} 
=
\mathcal J^{(-k)} 
\tilde{\longrightarrow} 
\mathcal J^{(-k-1)} = \Omega(\A, 2k+1)^{W}.
$$
\end{cor} 

\medskip

Concerning the strictly increasing filtration
\[
\dots
\Omega(\A, 2k-1) \subset 
\Omega(\A, 2k) \subset 
\Omega(\A, 2k+1) \subset 
\dots,
\]
the following Proposition asserts the $W$-invariant
parts of 
$\Omega(\A, 2k-1)$
and
$
\Omega(\A, 2k)
$  
are equal.

\begin{prop} 
$\Omega(\A, 2k)^{W} = \Omega(\A, 2k-1)^{W}=
\mathcal J^{(-k)}   $
for $k\in\Z$. 
In particular,
$
\Omega_{R} = 
\Omega_{S}^{W} = \Omega(\A, -1)^{W} . $ \end{prop} 

\medskip

\noindent
{\bf Proof.} 
It is obvious that
$\Omega(\A, 2k-1) \subseteq \Omega(\A, 2k)$
because 
$
R_{-2k+1} = R_{-2k} J(\P) (B^{(1-k)})^{-1} B
$ by Proposition \ref{inductive} (3).
Thus 
$\Omega(\A, 2k-1)^{W}  \subseteq \Omega(\A, 2k)^{W} $.

Let 
$
\omega = \sum_{j=1}^{\ell} f_{j} \,\omega^{(-2k)}_{j}
\in
\Omega(\A, 2k)^{W}$
with $f_{j} \in S$.  
Since
\begin{align*} 
({\rm
Eq})_{k}
\,\,\,\,\,
\,\,\,\,\,
\left[
\omega^{(-2k)}_{1}, 
\dots ,
\omega^{(-2k)}_{\ell} 
\right]
&=
\left[
\omega^{(-2k-1)}_{1}, 
\dots ,
\omega^{(-2k-1)}_{\ell} 
\right]
D[J(\P)]^{-1}
\end{align*}
by Proposition \ref{inductive} (2), we may express
\begin{align*} 
\omega 
= 
\sum_{j=1}^{\ell} f_{j} \,\omega^{(-2k)}_{j}
=
\sum_{j=1}^{\ell} f_{j} 
\left(
\sum_{i=1}^{\ell} h_{ij} \,
\omega_{i}^{(-2k-1)}
\right)
=
\sum_{i=1}^{\ell} \left(
\sum_{j=1}^{\ell} h_{ij} \,
f_{j} 
\right)
\omega_{i}^{(-2k-1)},    
\end{align*} 
where
$h_{ij} $ is the $(i, j)$-entry of $D[J(\mathbf P)]^{-1} $.
Note that
$\omega\in\Omega(\A, 2k+1)^{W} $ 
and
that
$\Omega(\A, 2k+1)^{W} $ 
has a basis
$
\{
\omega^{(-2k-1)}_{1},
$
$
\omega^{(-2k-1)}_{2},
\dots
,
\omega^{(-2k-1)}_{\ell}
\}   $ over $R$.
Then
we know that 
$
\sum_{j=1}^{\ell} 
h_{ij} \, f_{j}
$
is $W$-invariant
for
$1\le i \le \ell$. 
Applying 
$({\rm
Eq})_{0}
$
we have
\begin{align*}
\omega'
:=
\sum_{j=1}^{\ell} f_{j} \, dx_{j}
=
\sum_{j=1}^{\ell} f_{j} \,\omega^{(0)}_{j}
&=
\sum_{j=1}^{\ell} 
f_{j} 
\sum_{i=1}^{\ell} 
\,
h_{ij}  
\,\omega^{(-1)}_{i}\\
&=
\sum_{i=1}^{\ell} 
\left(
\sum_{j=1}^{\ell} 
h_{ij}  
\,
f_{j} 
\right)
\omega^{(-1)}_{i}
\in
\Omega_{S}^{W}.  
\end{align*}
Recall $\Omega_{S}^{W} = \Omega_{R} =
\oplus_{i=1}^{\ell} R\,\left( dP_{i} \right)  $ 
by \cite{sol1}. 
Thus there exist 
$g_{i} \in R\,\,\,
(1\le i\le \ell)$
 such that
$$
\omega' = 
\sum_{i=1}^{\ell} g_{i} \, \left(dP_{i}\right)
=
\sum_{j=1}^{\ell} 
\left(
\sum_{i=1}^{\ell} 
g_{i} \, \left(\partial P_{i} / \partial x_{j}\right) 
\right)
dx_{j}.
$$
 This implies
\[
f_{j} 
=
\sum_{i=1}^{\ell} 
g_{i} \, \left(\partial P_{i} / \partial x_{j}\right)
\,\,\,
(1\le i\le \ell). 
\]
%
Since
\[
\left[
\omega^{(-2k)}_{1}, 
\dots ,
\omega^{(-2k)}_{\ell} 
\right]
J(\P)
=
\left[
\omega^{(-2k+1)}_{1}, 
\dots ,
\omega^{(-2k+1)}_{\ell} 
\right]
B^{-1} B^{(1-k)}
\]
by Proposition \ref{inductive} (3),
one has
\begin{align*} 
\omega
&=
\sum_{j=1}^{\ell} 
f_{j} 
\,\omega^{(-2k)}_{j} 
=
\sum_{j=1}^{\ell} 
\left(
\sum_{i=1}^{\ell} 
g_{i} \, \left(\partial P_{i} / \partial x_{j}\right)
\right)
\,\omega^{(-2k)}_{j}\\ 
&=
\sum_{i=1}^{\ell}
 g_{i} 
\left(
\sum_{j=1}^{\ell}
(\partial P_{i}/\partial x_{j}) \,\omega^{(-2k)}_{j}
\right)
\in
\bigoplus_{i=1}^{\ell} 
R \,\, \omega^{(-2k+1)}_{i}
=
\Omega(\A, 2k-1)^{W}.   
\end{align*} 
This proves
$
\Omega(\A, 2k)^{W}
\subseteq
\Omega(\A, 2k-1)^{W}.
$ 
\owari

\section{The case of  derivations}
\label{section3} 
Denote 
$\partial/\partial x_{i} $ 
and
$\partial/\partial P_{i} $ 
simply by 
$\partial_{x_{i}}  $ 
and
$\partial_{P_{i} } $ 
respectively.
Then 
\[
\Der_{S} = \bigoplus_{j=1}^{\ell} S \,\partial_{x_{j}},
\ \ 
\Der_{R} = \bigoplus_{j=1}^{\ell} R \,\partial_{P_{j}},
  \ \ 
\Der_{F} = \bigoplus_{j=1}^{\ell} F \,\partial_{x_{j} }.
\]

In this section we translate the results in the 
previous section by the $F$-isomorphism
\[
I^{*} : \Omega_{F} \rightarrow \Der_{F} 
\]
defined by
$I^{*} (\omega) (f) = I^{*} (\omega, df)$ for
$f\in F$ and $\omega\in\Omega_{F}  $.
Explicitly we can express
\[
I^{*} 
(\sum_{j=1}^{\ell} f_{j} \ dx_{j} )
=
\sum_{j=1}^{\ell}
\left(
\sum_{i=1}^{\ell}  I^{*} (dx_{i}, dx_{j}) \,f_{i} \right) \partial_{x_{j}}
\]
for $f_{j} \in F \ \ (1\le j\le \ell)$.

\begin{define}
Define 
$
\eta_j^{(m)}
:=
I^{*}(\omega^{(m)}_{j})
$
for 
$m\in\Z, \ 1\le j\le\ell$. 
\label{Qeta}
\end{define}

Then
\[
[ \eta_1^{(m)}, \ldots,\eta_\ell^{(m)} ]
=
[ \partial_{x_1},\ldots, \partial_{x_\ell} ] A R_{m}. 
\]
In particular, 
\begin{eqnarray*}
[ \eta_1^{(1)}, \ldots,\eta_\ell^{(1)} ]
=
[ \partial_{x_1},\ldots, \partial_{x_\ell} ] A J(\mathbf{P})
=
[I^*(dP_1),\ldots,I^*(dP_\ell)],
\end{eqnarray*}
\begin{eqnarray*}
[ \eta_1^{(-1)}, \ldots,\eta_\ell^{(-1)} ] 
&=&
[ \partial_{x_1},\ldots, \partial_{x_\ell} ] AD[J(\mathbf{P})]
=
[ \partial_{x_1},\ldots, \partial_{x_\ell} ] J(\mathbf{P})^{-T} B \\ 
&=&
[ \partial_{P_1},\ldots,\partial_{P_\ell}] B.
\end{eqnarray*}

\begin{define}
Define 
\begin{align*}
D(\A, m):= 
\{
\theta\in \Der_{S} 
\mid
\theta(\alpha_{H})
\in
S \cdot \alpha_{H}^{m} 
\  \text{~for all~}
H\in \A
 \}
\end{align*}
for $m\ge 0$ which is the $S$-module of {\bf
logarithmic derivations}
along $\A$ of contact order $m$. 
When $m < 0$ define
\[
D(\A, m)
:=
\bigoplus_{1\le j\le\ell} S\,\eta_{j}^{(m)}.
\]
Lastly define
\[
D(\A, -\infty)
:=
\bigcup_{m\in\Z}  D(\A, m).
\]

\label{multiderivation}
\end{define}

\begin{theorem} 
$D(\A, m)$ is a free $S$-module
with a basis 
$
\eta^{(m)}_{1},
\eta^{(m)}_{2},
\dots,
\eta^{(m)}_{\ell}
$ for $m\in\Z.$   
\label{SbasisD}
\end{theorem} 

\medskip

\noindent
{\bf Proof.}
 {\it Case 1.} When $m<0$ this is nothing but the definition.

 {\it Case 2.} Let $m\ge 0$.
For a canonical contraction $\langle\ ,\ \rangle: 
\Der_F \times \Omega_F \rightarrow F$, 
define
the $(\ell \times \ell)$-matrix 
$$
Y_{m} 
:=
[ \langle \omega_i^{(-m)}, \eta_j^{(m)} \rangle ]_{1 \le i,j \le \ell} 
=
R_{-m} A R_{m} 
$$
for
$ m \ge 0$. 
Since
the two $S$-modules $\Omega(\A, m)$ and $D(\A, m)$ are
dual each other (see \cite{Z}) , it is enough to show
that $\det Y_{m} \in \mbox{GL}_{\ell} (S)$. 
It follows from the following Proposition \ref{Ym}.
\owari

\begin{cor} 
$I^{*} (\Omega(\A, m)) = D(\A, -m)$ for $m\in\Z$
and
$I^{*} (\Omega(\A, \infty)) = D(\A, -\infty)$. 
\label{do1}
\end{cor}

\begin{cor} 
$\Omega(\A, -m) = 
\{
\omega\in\Omega_{S} 
\mid
I^{*} (\omega, d\alpha_{H} )
\in S \cdot \alpha_{H}^{m} 
\text{~for any~ } H\in\A \}
$
for $m>0$.
\label{do2}
\end{cor}

\begin{prop} 
\label{Ym} 
\hspace{5mm} (1) $Y_{2k-1} = (-1)^{k+1} 
B^{T} (B^{(k)})^{-1} B \in \mbox{GL}_{\ell} (T)$ for $k\in\Z$, 

(2) $Y_{2k} = (-1)^{k} A \in \mbox{GL}_{\ell} (\R)$ for $k\in\Z$.
\end{prop} 

\medskip

\noindent
{\bf Proof.} 

(1)
{\it Case 1.1.}
Let $m=2k-1$ with $k\ge 1$. 
We prove by an induction on $k$. 
When $k=1$, 
\begin{eqnarray*}
Y_{1} = R_{-1}^T A R_1= 
D[J(\mathbf{P})]^T A J(\mathbf{P})
&=&B^T\in \mbox{GL}_{\ell} (T).
\end{eqnarray*}
Assume that $k>1$ and prove by induction. By using Proposition \ref{inductive}
(5) and (4), we obtain 
\begin{eqnarray*}
Y_{2k-1} 
&=&R_{1-2k}^TAR_{2k-1}= D[R_{3-2k}]^T AR_{2k-3}B^{-1}G(B^{(k)})^{-1} B\\
&=& \{D[R_{3-2k}^T AR_{2k-3} ] -R_{3-2k}^T D[AR_{2k-3}] \} B^{-1}G(B^{(k)})^{-1}B \\
&=& -R_{3-2k}^T AR_{2k-5} B^{-1} G(B^{(k-1)})^{-1} B B^{-1} B^{(k-1)}(B^{(k)})^{-1} B\\
&=& -R_{3-2k}^T AR_{2k-3} B^{-1} B^{(k-1)}(B^{(k)})^{-1} B\\
&=&(-1)^{k+1} B^T(B^{(k-1)})^{-1} B B^{-1} B^{(k-1)} (B^{(k)} )^{-1} B\\
&=&(-1)^{k+1} B^T (B^{(k)} )^{-1}B.
\end{eqnarray*}

{\it Case 1.2.} 
Next assume that $m=2k-1$
with $ k \le 0$. Recall that 
$$
(B^{(1-k)})^T= -k B+(1-k)B^T=-B^{(k)}.
$$
Then
\begin{eqnarray*}
R^T_{1-2k}AR_{2k-1}
&=& (R_{2k-1}^T AR_{1-2k})^T 
= ((-1)^{k}B^T (B^{(1-k)})^{-1}  B)^{T}  \\
&=& (-1)^{k+1}B^T(B^{(k)})^{-1} B.
\end{eqnarray*}

(2)
Apply (1),
Proposition \ref{inductive} (2) and (3) to compute
\begin{align*} 
R^T_{-2k}AR_{2k}
&=
J(\P)^{-T} (B^{(1-k)})^{T} B^{-T} 
R_{-2k+1}^{T} A R_{2k-1} 
B^{-1} J(\P)^{T} A\\
&=
J(\P)^{-T} (B^{(1-k)})^{T} B^{-T} 
Y_{2k-1} 
B^{-1} J(\P)^{T} A
=
(-1)^{k} A.
\ \ \ 
\square
\end{align*}

\noindent
{\it Remark.} 
Corollaries \ref{do1} and \ref{do2} show that 
the definitions of $D(\A,m)$ and $\Omega(\A,m)$ for $m \in \Z_{<0}$ are 
equivalent to those of $D\Omega(\A,m)$ and $\Omega D(\A,m)$ in \cite{A4}. 
\medskip

Consider the $T$-linear connection (covariant derivative)
$$
\nabla_{D} :
\Der_{F} 
\rightarrow
\Der_{F} 
$$  
characterized by
$\nabla_{D} (f X)
=
(Df) X + f (\nabla_{D} X)$ 
and 
$\nabla_{D} (\partial_{x_{j} }  ) =
0$ 
for $f\in F$,
$X\in\Der_{F} $ 
and $1\le j\le \ell$.  
Then it is easy to see the diagram
$$\xymatrix@R1pc{
  \Omega_{F}  \ar[r]^{\nabla_D} \ar[d]_{I^{*}} 
&\Omega_{F} \ar[d]_{I^{*}}\\
  \Der_{F}     \ar[r]^{\nabla_D} 
&\Der_{F}}     
$$
is commutative. In fact
\begin{align*} 
\nabla_{D}\circ I^{*} 
\left(
\sum_{j=1}^{\ell} 
f_{j} \, dx_{j}  
\right) 
&=
\nabla_{D}
\left[
\sum_{j=1}^{\ell} 
\left(
\sum_{i=1}^{\ell} 
I^{*}(dx_{i}, dx_{j})
f_{i}   
\right)
\, \partial_{x_{j}}
\right] \\
&=
\sum_{j=1}^{\ell} 
\left(
\sum_{i=1}^{\ell} 
I^{*}(dx_{i}, dx_{j})
D(f_{i})   
\right)
\, \partial_{x_{j}}\\
&=
I^{*} 
\left(
\sum_{j=1}^{\ell} 
D(f_{j})   
\, dx_{j}
\right)
=
I^{*} 
\circ 
\nabla_{D}
\left(
\sum_{j=1}^{\ell} 
f_{j} \, dx_{j}  
\right) .
\end{align*}

Define
$
\mathcal C_{k} :=
I^* (\mathcal B_{k-1} )
=
\{
\eta^{(2k-1)}_{1},
\eta^{(2k-1)}_{2},
\dots
,
\eta^{(2k-1)}_{\ell}\}
$ 
for each
$k\in\Z$.
The following Theorems \ref{3.4} and
\ref{3.1} can be proved by translating Theorems \ref{0.4} 
and \ref{0.1} through
$\nabla_{D} $.

\begin{theorem}
\label{3.4}

\hspace{5mm} 
(1)
The $R$-module 
$D(\A, 2k-1)^{W} $
is free
with a basis
$\mathcal C_{k} $ 
for $k\in\Z$.

(2)
The $T$-module 
$D(\A, 2k-1)^{W} $
is free
with a basis
$\bigcup_{p\ge k} \mathcal C_{p} $ 
for $k\in\Z$.

(3)
$
\mathcal C
:=
\bigcup_{k\in\Z} \mathcal C_{k} 
$ 
is a basis for 
$D(\A, -\infty)^{W} $
as a $T$-module.
\end{theorem}

\begin{define}
Define
$$
\mathcal G_{k} 
:=
I^{*}(\mathcal F_{k-1}),
\ \
\mathcal H^{(k)} 
:=
I^{*}(\mathcal J^{(k-1)})\ \ (k\in\Z, \ 1\le j\le\ell).
$$
\end{define}

Then
\[
\mathcal G_{k} = \bigoplus_{1\le j\le \ell} T\, 
\eta_{j}^{(2k-1)},
\ \
\mathcal H^{(k)} 
= \bigoplus_{p\ge k} \,\mathcal G_{p}.
\]
The $\nabla_{D} $ induces $T$-isomorphisms
\[
\nabla_{D} : \mathcal G_{k+1} \tilde{\longrightarrow} \mathcal G_{k},
\ \ \
\nabla_{D} : D(\A, 2k+1)^W \tilde{\longrightarrow} D(\A, 2k-1)^W.
\]

In particular,
\[
\mathcal G_{0} = \bigoplus_{j=1}^{\ell} T \, \partial_{P_{j} },
\ \ \ \ \text{and} \ \ \ \ 
\mathcal H^{(0)} = \bigoplus_{j=1}^{\ell} R \, 
\partial_{P_{j} } = \Der_{R}.
\]

\begin{theorem}
\label{3.1}

\hspace{5mm} 
(1)
The $\nabla_{D} $ induces a $T$-linear automorphism  
$\nabla_{D} : 
D(\A, -\infty)^{W} 
\stackrel{\sim}{\rightarrow} 
D(\A, -\infty)^{W}. 
$ 

(2)
$
D(\A, -\infty)^{W} 
=
{
\bigoplus_{k\in\Z}}\, {\mathcal G}_{k}.
$

(3)
$ D(\A, 2k-1)^{W}
=
\mathcal H^{(k)} =
\bigoplus_{p\ge k}\, {\mathcal G}_{p}.
\,\,\, (k \in\Z)$.
\end{theorem}

\noindent
{\it Remark.}
The construction of 
a basis
$
\eta^{(1)}_{1}, \dots
,
\eta^{(1)}_{\ell}$ 
for
$D(\A, 1)$ 
is
due to K. Saito \cite{Sa1}. 
A basis for 
$D(\A, 2)$
was constructed in
\cite{ST2}. 
In
\cite{T4}
$D(\A, m)$ was found to be 
a free 
$S$-module for all $m\ge 0$ whenever
$\A$ is a Coxeter arrangement. 
Note that it is re-proved in Theorem 
\ref{SbasisD} in this article.
In 
\cite{Sa4} 
K. Saito called the decreasing filtration $
\Der_{R} 
=
\mathcal H^{(0)} 
\supset
\mathcal H^{(1)} 
\supset
\dots
$
and the decomposition
$\Der_{R} = D(\A, -1)^{W}
= 
\mathcal H^{(0)} =
\bigoplus_{p\ge 0}\, {\mathcal G}_{p}
$ 
the Hodge filtration
and 
the Hodge decomposition respectively.
They are essential to define the flat structure
(or equivalently the Frobenius manifold structure in 
topological field theory) on the orbit space $V/W$.
Note that
Theorem \ref{3.1} (3), when $k\ge 0$,
is the main theorem of \cite{T6}.

\section{Relation among bases for logarithmic forms and derivations}
\label{section4}
 In the previous section we constructed a basis
$\{\omega^{(m)}_{j}\}$ for $\Omega(\A, m)$ and a
basis 
$\{\eta^{(m)}_{j}\}$
for
$D(\A,m)$ for 
$m\in\Z$. 
In this section we briefly describe their relations to other 
bases
constructed in the
earlier works
\cite{T4}, 
\cite{Y0}, and 
\cite{AY2}. 
In 
\cite{T4},
the following bases for $D(\A, 2k+1)$ 
and $D(\A, 2k)$ are given:
\begin{align*} 
[\xi_1^{(2k+1)}, \ldots,\xi_\ell^{(2k+1)}]
&:=[\partial_{x_1},\ldots,\partial_{x_\ell} ] 
A J(D^k[\mathbf{x}])^{-1} J(\mathbf{P}),\\
[\xi_1^{(2k)}, \ldots, \xi_\ell^{(2k)}]
&:=
[\partial_{x_1},\ldots,
\partial_{x_\ell} ] A J(D^k[\mathbf{x}])^{-1}. 
\end{align*} 
The
two bases $\{\eta_{j}^{(m)}\}$ and $\{\xi_{j}^{(m)}\}$ are 
related as follows:

\begin{prop}
For $k \in \Z_{\ge 0}$, 
\begin{align*} 
[\xi_1^{(2k+1)},\ldots,\xi_\ell^{(2k+1)} ] &=
(-1)^k [\eta_1^{(2k+1)},\ldots,\eta_\ell^{(2k+1)} ] B^{-1} B^{(k+1)},\\
[\xi_1^{(2k)},\ldots,\xi_\ell^{(2k)} ] &=
(-1)^k [\eta_1^{(2k)},\ldots,\eta_\ell^{(2k)} ].
\end{align*} 
\label{relation1}
\end{prop}

\noindent
\textbf{Proof}. 
The second formula is immediate from Definition
\ref{matrixR}. 
The following computation proves 
 the first formula:
\begin{eqnarray*}
J(D^k[\mathbf{x}])^{-1}J(\mathbf{P})&=&
(-1)^{k+1} R_{2k+1} D[J(\mathbf{P})]^{-1} J(D^{k+1}[\mathbf{x}]) 
J(D^k[\mathbf{x}])^{-1} J(\mathbf{P})\\
&=& 
(-1)^{k} R_{2k+1} D[J(\mathbf{P})]^{-1} A^{-1} J(\mathbf{P})^{-T} B^{(k+1)} \\
&=& 
(-1)^{k} R_{2k+1} B^{-1} B^{(k+1)}.
\ \ \square
\end{eqnarray*}

In 
\cite{Y0},
the following bases 
are given:
\begin{align*} 
[\nabla_{I^{*}(dP_{1})} \nabla_D^{-k}  \theta_E, \ldots,
\nabla_{I^{*}(dP_{\ell})} \nabla_D^{-k}  \theta_E]
\ \ \ &\text{for}\ \ D(\A, 2k+1),\,\\
[\nabla_{\partial_{x_{1}}} \nabla_D^{-k}  \theta_E, \ldots,
\nabla_{\partial_{x_{\ell}}} \nabla_D^{-k}  \theta_E]
\ \ \ &\text{for}\ \ D(\A, 2k).
\end{align*} 
Here $\theta_E$ is the Euler derivation.
Their relations to
$\{\eta_{j}^{(m)}\}$ are 
given as follows:

\begin{prop}
Let $k \in \Z_{\ge 0}$. Then 
\begin{align*} 
[\nabla_{I^*(dP_1)} \nabla_D^{-k}  \theta_E, \ldots,
\nabla_{I^*(dP_\ell)} \nabla_D^{-k}  \theta_E]
&=
[\eta_1^{(2k+1)},\ldots,\eta_\ell^{(2k+1)}] B^{-1} B^{(k+1)},\\
[\nabla_{\partial_{x_{1}}} \nabla_D^{-k}  \theta_E, \ldots,
\nabla_{\partial_{x_{\ell}}} \nabla_D^{-k}  \theta_E]
&=[\eta_1^{(2k)},\ldots,\eta_\ell^{(2k)}] A^{-1} .
\end{align*} 
\label{deri}
\end{prop}

\noindent
{\bf Proof.}
By \cite[Theorem 1.2.]{T5} and 
\cite{T7} one has
$$
[\nabla_{I^*(dP_1)} \nabla_D^{-k}  \theta_E, \ldots,
\nabla_{I^*(dP_\ell)} \nabla_D^{-k}  \theta_E]
=
(-1)^k
[\xi_1^{(2k+1)},\ldots,\xi_\ell^{(2k+1)}].
$$
Combining with Proposition \ref{relation1},
we have the first relation.
For the second one, compute
\begin{align*} 
[\nabla_{\partial_{x_{1}}} \nabla_D^{-k}  \theta_E, \ldots,
\nabla_{\partial_{x_{\ell}}} \nabla_D^{-k}  \theta_E]A J(\P)
&=
[\nabla_{I^*(dP_1)} \nabla_D^{-k}  \theta_E, \ldots,
\nabla_{I^*(dP_\ell)} \nabla_D^{-k}  \theta_E]\\
&
=
[\eta_1^{(2k+1)},\ldots,\eta_\ell^{(2k+1)}] B^{-1} B^{(k+1)}\\
&
=
[\eta_1^{(2k)},\ldots,\eta_\ell^{(2k)}] 
J(\P)
\end{align*} 
by Proposition \ref{inductive} (3).
\owari

\medskip

Next let us review 
the bases for $\Omega(\A, m)$ 
described in 
\cite[Theorem 6]{AY2}: 
Let $k \in \Z_{\ge 0}$ and 
$P_1$ the smallest degree 
basic invariant.
Then
$$
\{\nabla_{\partial_{P_1}} \nabla_D^{k}  dP_1, \ldots,
\nabla_{\partial_{P_\ell}} \nabla_D^{k}  dP_1\}
$$
forms a basis for $\Omega(\A,2k+1)$
and 
$$
\{\nabla_{\partial_{x_1}} \nabla_D^{k}  dP_1, \ldots,
\nabla_{\partial_{x_\ell}} \nabla_D^{k}  dP_1\}
$$
forms a basis for $\Omega(\A,2k)$.

\begin{prop}
Let $k \ge 0$. Then 
\begin{align*} 
[\nabla_{\partial_{P_1}} \nabla_D^{k}  dP_1, \ldots,
\nabla_{\partial_{P_\ell}} \nabla_D^{k}  dP_1]
&=
[\omega_1^{(-2k-1)},\ldots,\omega_\ell^{(-2k-1)}]
B^{-1}, \\
[\nabla_{\partial_{x_1}} \nabla_D^{k}  dP_1, \ldots,
\nabla_{\partial_{x_\ell}} \nabla_D^{k}  dP_1]
&
=
[\omega_1^{(-2k)},\ldots,\omega_\ell^{(-2k)}]A^{-1}. 
\end{align*}

\label{commu}
\end{prop}

\noindent
\textbf{Proof}. 
First, note that $[\nabla_D,\nabla_{\partial_{P_i}} ]$ is $W$-invariant, 
hence in $\Der_R$. Since the smallest degree of derivations in $\Der_R$ is 
$\deg \partial_{P_\ell}$, it follows that $[\nabla_D,\nabla_{\partial_{P_i}}]=0$. 
In other words, 
$\nabla_{\partial_{P_i}}$ and 
$\nabla_{\partial_{P_\ell}}=\nabla_D$ commute for all $i$. Hence 
$$
[\nabla_{\partial_{P_1}} \nabla_D^{k}  dP_1, \ldots,
\nabla_{\partial_{P_\ell}} \nabla_D^{k}  dP_1]
=\nabla_D^k[\nabla_{\partial_{P_1}}  dP_1, \ldots,
\nabla_{\partial_{P_\ell}}   dP_1].
$$
Our proof is an induction on $k$. 
First assume that $k=0$. Choose 
$$
P_1=\displaystyle \frac{1}{2}[x_1,\ldots,x_\ell]A^{-1}[x_1,\ldots,x_\ell]^T,
$$
and
$$
dP_{1} = [dx_1,\ldots,dx_\ell]A^{-1}[x_1,\ldots,x_\ell]^T.
$$
Compute
\begin{eqnarray*}
[\nabla_{\partial_{P_1}}  dP_1, \ldots,
\nabla_{\partial_{P_\ell}}   dP_1]B
&=&
[\nabla_{\partial_{x_1}}  dP_1, \ldots,
\nabla_{\partial_{x_\ell}}   dP_1]J(\P)^{-T} B
\\
&=&[dx_1,\ldots,dx_\ell]A^{-1}J(\mathbf{P})^{-T} B\\
&=&[dx_1,\ldots,dx_\ell]D[J(\mathbf{P})]
=
[\omega_1^{(-1)},\ldots,\omega_\ell^{(-1)}].
\end{eqnarray*}
For $k>0$, apply $\nabla_D^k$ and use the commutativity. 
Then we have the first relation.
For the second relation use Proposition \ref{inductive} (2)
to compute:
\begin{eqnarray*}
[\nabla_{\partial_{x_1}}  \nabla_D^k dP_1, \ldots,
\nabla_{\partial_{x_\ell}}  \nabla_D^k  dP_1]
&=&
[\nabla_{\partial_{P_1}} \nabla_D^k dP_1, \ldots,
\nabla_{\partial_{P_\ell}} \nabla_D^k dP_1] 
J(\P)^{T}\\ 
&=&
[\omega_1^{(-2k-1)},\ldots,\omega_\ell^{(-2k-1)}]
B^{-1} 
J(\P)^{T}\\ 
&=&[dx_1,\ldots,dx_\ell]
R_{-2k-1} 
B^{-1} 
J(\P)^{T}\\ 
&=&[dx_1,\ldots,dx_\ell]R_{-2k} A^{-1} \\
&=&
[\omega_1^{(-2k)},\ldots,\omega_\ell^{(-2k)}]A^{-1}.
\end{eqnarray*}
\owari
\medskip

\noindent
{\it Remark.} 
If $k <0$ in Propositions \ref{deri} and \ref{commu}, then 
the derivations and $1$-forms in the left hand sides are proved 
to form 
bases for the logarithmic modules $D\Omega(\A,2k+1),D\Omega(\A,2k), 
\Omega D(\A,2k+1)$ and $\Omega D(\A,2k)$ in \cite{A4}. By using the 
same arguments in the proofs above, we can show that 
Propositions \ref{deri} and \ref{commu} hold true 
for all integers $k$ in the logarithmic modules $D\Omega(\A,\textbf{m})$ 
and $\Omega D(\A,\textbf{m})$ with $\textbf{m}:\A \rightarrow \Z$.


%
%
%
%
%
%
%
%

 \vspace{5mm}

\end{document}